\documentclass[12pt]{article}

\usepackage{hyperref}
\usepackage{amsmath}
\usepackage{amssymb}
\usepackage{amsfonts}
\usepackage{amsopn}
\usepackage{amsthm}
\usepackage[all]{xy}

\swapnumbers
\theoremstyle{plain}
\newtheorem{thm}{Theorem}[section]

\newtheorem{lem}[thm]{Lemma}

\theoremstyle{definition}
\newtheorem{ntt}[thm]{}

\newtheorem{rem}[thm]{Remark}
\newtheorem{dfn}[thm]{Definition}

\newcommand{\pl}{\mathbb{P}}   
\newcommand{\zz}{\mathbb{Z}}   
\newcommand{\A}{\mathrm{A}} 
\newcommand{\id}{\mathrm{id}}       
\newcommand{\pr}{\mathrm{pr}}       
\newcommand{\op}{\mathrm{op}}     

\newcommand{\Var}{\mathcal{V}ar} 
\newcommand{\M}{\mathcal{M}}     
\newcommand{\E}{\mathcal{E}} 
\newcommand{\Cor}{\mathcal{C}or} 
\newcommand{\zab}{\mathbb{Z}\text{-}\mathcal{A}b} 

\newcommand{\xra}[1]{\xrightarrow{#1}} 

\DeclareMathOperator{\Br}{\mathrm{Br}}     
\DeclareMathOperator{\CH}{\mathrm{CH}}      
\DeclareMathOperator{\Mor}{\mathrm{Mor}}    
\DeclareMathOperator{\im}{\mathrm{im}}      
\DeclareMathOperator{\SB}{\mathrm{SB}}      
\DeclareMathOperator{\Gr}{\mathbb{G}} 
\DeclareMathOperator{\ind}{\mathrm{ind}} 

\title{Motivic decomposition of a generalized Severi-Brauer variety}

\author{K.Zainoulline}

\date{\today}

\begin{document}

\maketitle

\begin{abstract}
Let $A$ and $B$ be two central simple algebras of a prime degree $n$ over
a field $F$
generating the same subgroup in the Brauer group $\Br(F)$.
We show that the Chow motive of a Severi-Brauer variety $\SB(A)$
is a direct summand of the motive of a generalized Severi-Brauer
variety $\SB_d(B)$
if and only if
$[A]=\pm d [B]$ in $\Br(F)$. The proof uses methods of Schubert calculus
and combinatorial properties
of Young tableaux, e.g., the Robinson-Schensted correspondence. 

\vspace{1ex}

\end{abstract}

\section{Introduction}

Let $X$ be a twisted flag $G$-variety for a linear algebraic group $G$.
The main result of paper \cite{CPSZ} says that under certain
restrictions the Chow motive of $X$ can be expressed in terms 
of motives of ``minimal'' flag varieties, i.e., 
those which correspond
to maximal parabolic subgroups of $G$. 
A natural question arises:
is it possible to decompose the motive of such a ``minimal'' flag?

A particular case of such a decomposition
was already provided in \cite{CPSZ}.
More precisely, a ``minimal flag'' $X$  was
a generalized Severi-Brauer variety of 
ideals of reduced dimension $2$ in a division algebra $B$ of degree $5$ 
and the decomposition 
was (see \cite[Theorem~2.5]{CPSZ}) 
$$
\M(\SB_2(B)) \simeq \M(\SB(A))\oplus\M(\SB(A))(2),
$$
where $[B]=2[A]$ in the Brauer group.

In the present paper we provide an affirmative answer on this question
for any adjoint group $G$ of inner type $\A_n$ of a prime rank $n$. 
Namely, we show 
that the motive of a generalized Severi-Brauer variety
always contains (as a direct summand) the motive of a Severi-Brauer variety.

\begin{thm}\label{introthm}
Let $A$ and $B$ be two central simple algebras of a prime degree $n$
over a field $F$ generating the same subgroup in the Brauer group $\Br(F)$.
Then the motive of a Severi-Brauer variety $\SB(A)$ 
is a direct summand of the motive of a 
generalized Severi-Brauer variety $\SB_d(B)$
if and only if
\begin{equation}\label{mcond}
[A]=\pm d\,[B] \text{ in } \Br(F).
\end{equation} 
\end{thm}

\begin{rem} Theorem~\ref{introthm} may also be considered 
as a generalization of the
result by N.~Karpenko (see \cite[Criterion~7.1]{Ka00}) which says that 
the motives of Severi-Brauer
varieties $\SB(A)$ and $\SB(B)$ of central simple algebras $A$ and $B$ 
are isomorphic if and only if $[A]=\pm [B]$ in $\Br(F)$. 
\end{rem}

\begin{rem} We expect that Theorem~\ref{introthm} holds
in the case when $d$ and $n$ are coprime (see Section~\ref{scop}). 
Observe that if $d$ and $n$ are not coprime, 
then the theorem fails. 
It can be already seen on the level of generating functions.
Namely, consider the Grassmann variety $\Gr_2(4)$ 
of $2$-planes in a $4$-dimensional affine space ($n=4$ and $d=2$). 
Then for the generating functions we have
$$
P(\Gr_2(4),t)=(t^2+1)(t^2+t+1)\text{ and }P(\pl^3,t)=t^3+t^2+t+1.
$$ 
Obviously $P(\pl^3,t)$ doesn't divide $P(\Gr_2(4),t)$.
Indeed, it can be shown that $P(\pl^{n-1},t)$ divides 
$P(\Gr_d(n),t)$ if and only if $d$ and $n$ are coprime.
\end{rem}

\begin{rem} For motives with $\zz/n\zz$-coefficients 
there is the following decomposition (see \cite[Proposition~2.4]{CPSZ})
\begin{equation}\label{motdec}
\M(\SB_d(B)) \simeq \bigoplus_i \M(\SB(B))(i)^{\oplus a_i},
\end{equation}
where the integers $a_i$ are coefficients of the quotient of 
Poincar\'e polynomials 
$\tfrac{P(\Gr_d(n),t)}{P(\pl^{n-1},t)}=\sum_i a_it^i$ (see \ref{poin}).
Note that in this case the motives
of Severi-Brauer varieties corresponding to different classes of algebras
generating the same subgroup in the Brauer group are isomorphic
(see \cite[Section~7]{Ka00}), i.e., $\M(\SB(A))\simeq\M(\SB(B))$.
\end{rem}

The proofs are based on Rost Nilpotence Theorem for
projective homogeneous varieties proved 
by  V.~Chernousov, S.~Gille and A.~Merkurjev in \cite{CGM}.
Briefly speaking, this result reduces the problem of decomposing
the motive of a variety $X$ over $F$ into the question
about algebraic cycles in the Chow ring $\CH(X_s\times X_s)$
over the separable closure $F_s$.
Namely, the motive of $\SB(A)$ is a direct summand of the
motive of $\SB_d(B)$
if there exist
two cycles $f$ and $g$ in $\CH(\pl^{n-1}\times\Gr_d(n))$ 
such that both cycles belong to the image of the restriction map
$\CH(\SB(A)\times\SB_d(B))\to\CH(\pl^{n-1}\times\Gr_d(n))$, 
and the correspondence product $g^t\circ f$ is the identity.

We define $f$ and $g$ to be the Schur functions of
total Chern classes of certain bundles on $\pl^{n-1}\times\Gr_d(n)$.
Using the language and properties of Schur functions
we show that the identity $g^t \circ f=\id$ is a direct consequence
of the Robinson-Schensted correspondence, one of classical combinatorial
facts about Young tableaux.

The paper is organized as follows.
In section~\ref{prelim} we remind several definitions and notation 
used in the proofs. 
These include Chow motives, rational cycles, and generalized
Severi-Brauer varieties. In section~\ref{subrat} we describe the
subgroup of rational cycles of the Chow group of the product
of two generalized Severi-Brauer varieties. Indeed, we provide an explicit
set of generators for this subgroup modulo $n$ in terms of Schur functions.
In section~\ref{sevbrauer} we use this description
for proving some known results on motives of Severi-Brauer varieties.
Section~\ref{gsevbrauer} is devoted to the proof of the main theorem.
In section~\ref{proofid} we prove the crucial congruence
used in the proof of the main theorem.
In the last section we discuss the case of $\Gr_2(n)$,
where $n$ is an odd integer (not necessarily prime).

\section{Preliminaries}\label{prelim}

In the present section 
we remind definition of the category of Chow motives over a field $F$ 
following \cite{Ma68} and \cite{Ka01}. We recall the notion of 
a rational cycle and state the Rost Nilpotence
Theorem for idempotents following \cite{CGM}.
We recall several auxiliary facts concerning 
generalized Severi-Brauer varieties 
following \cite{KMRT98}.

\begin{ntt}[Chow motives]
Let $F$ be a field and $\Var_F$ be the category
of smooth projective varieties over $F$.
We define the category $\Cor_F$ of \emph{correspondences} over $F$.
Its objects are non-singular projective varieties over $F$.
For morphisms, called correspondences,
we set $\Mor(X,Y):=\CH^{\dim X}(X\times Y)$.
For two correspondences $\alpha\in \CH(X\times Y)$ and
$\beta\in \CH(Y\times Z)$ we define the composition
$\beta\circ\alpha\in \CH(X\times Z)$
$$
\beta\circ\alpha ={\pr_{13}}_*(\pr_{12}^*(\alpha)\cdot \pr_{23}^*(\beta)),
$$
where $\pr_{ij}$ denotes the projection
on product of the $i$-th and $j$-th factors of $X\times Y\times Z$
respectively and ${\pr_{ij}}_*$, ${\pr_{ij}^*}$ denote
the induced push-forwards and pull-backs for Chow groups.
The composition $\circ$
induces a ring structure on the abelian group $\CH^{\dim X}(X\times X)$.
The unit element of this ring is the class of diagonal cycle $\Delta_X$.

The pseudo-abelian completion of $\Cor_F$ is called the category
of \emph{Chow motives} and is denoted by $\M_F$.
The objects of $\M_F$
are pairs $(X,p)$, where $X$ is a non-singular projective variety
and $p$ is a projector, that is, $p\circ p=p$.
The motive $(X,\Delta_X)$ will be denoted by $\M(X)$.

By construction $\M_F$ is a self-dual tensor additive category,
where the duality is given by the transposition of cycles $\alpha\mapsto \alpha^t$
and the tensor product is given by the usual fiber product
$(X,p)\otimes(Y,q)=(X\times Y, p\times q)$.
Moreover, the contravariant Chow functor
$\CH: \Var_F \to \zab$ (to the category of $\zz$-graded abelian groups)
factors through $\M_F$, i.e., one has the commutative diagram of functors
$$
\xymatrix{
\Var_F \ar[rr]^{\CH}\ar[rd]_{\Gamma}& & \zab \\
 & \M_F \ar[ru]_{R}&
}
$$
where $\Gamma: f\mapsto \Gamma_f$ is the (contravariant) graph functor 
and $R:\M_F \to \zab$ is the (covariant) realization functor
given by $R:(X,p) \mapsto \im(p^*)$,
where $p^*$ is the composition
$$
p^*: \CH(X) \xra{\pr_1^*} \CH(X\times X) \xra{\cdot p} 
\CH(X\times X) \xra{{\pr_2}_*} \CH(X).
$$

 Consider the morphism
$(\id,e):\pl^1\times\{pt\}\to \pl^1\times\pl^1$. The image of
the induced push-forward $(\id,e)_*$ doesn't depend on the choice of a point
$e:\{pt\}\to \pl^1$
and defines the projector in $\CH^1(\pl^1\times\pl^1)$ denoted by $p_1$.
The motive $L=(\pl^1,p_1)$ is called {\em Lefschetz motive}.
For a motive $M$ and a nonnegative integer $i$ we denote by $M(i)=M\otimes L^{\otimes i}$ its {\em twist}.
Observe that
$$
\Mor((X,p)(i),(Y,q)(j))=q\circ \CH^{\dim X+i-j}(X\times Y)\circ p.
$$
\end{ntt}

\begin{ntt}[Product of cellular varieties] 
Let $G$ be a split linear algebraic group over a field $F$.
Let $X$ be a projective $G$-homogeneous variety, i.e., 
$X=G/P$, where $P$ is a parabolic subgroup of $G$.
The abelian group structure of $\CH(X)$, as well as
its ring structure, is well-known.
Namely, $X$ has a cellular filtration and the generators 
of Chow groups of the bases of this filtration
correspond to the free additive generators of $\CH(X)$ (see \cite{Ka01}).
Note that the product of two projective homogeneous varieties
$X\times Y$ has a cellular filtration as well,
and $\CH^*(X\times Y)\cong \CH^*(X)\otimes \CH^*(Y)$
as graded rings.
The correspondence product of two cycles 
$\alpha=f_\alpha \times g_\alpha \in \CH(X\times Y)$ and
$\beta=f_\beta \times g_\beta \in \CH(Y\times X)$ is given
by
\begin{equation}\label{composit}
(f_\beta\times g_\beta)\circ(f_\alpha\times g_\alpha)=
\deg(g_\alpha \cdot f_\beta)(f_\alpha\times g_\beta),
\end{equation}
where $\deg: \CH(Y)\to \CH(\{pt\})=\zz$ is the degree map.
\end{ntt}

\begin{ntt}[Rational cycles]
Let $X$ be a projective variety of dimension $n$ over $F$.
Let $F_s$ denote the separable closure of $F$.
Consider the scalar extension $X_s=X\times_F F_s$.
We say a cycle $J\in \CH(X_s)$ is {\it rational}
if it lies in the image of the pull-back homomorphism
$\CH(X)\to \CH(X_s)$.
For instance, there is an obvious rational cycle $\Delta_{X_s}$ on
$\CH^n(X_s\times X_s)$ that is given by the diagonal class.

Let $E$ be a vector bundle over $X$. Then the
total Chern class $c(E_s)$ of the pull-back induced by
the scalar extension $F_s/F$ is rational.
Let $L/F$ be a finite separable field extension which splits $X$, 
i.e., there is an induced isomorphism $\CH(X_L)\xra{\simeq} \CH(X_s)$.
Then the cycle $\deg(L/F)\cdot J$, $J\in \CH(X_s)$, is rational.
Observe that all linear combinations,
intersections and correspondence products of rational cycles
are rational.
\end{ntt}

\begin{ntt}[Rost nilpotence]\label{exproj}
We will use the following fact (see \cite[Cor.~8.3]{CGM})
that follows from the Rost Nilpotence Theorem.
Let $X$ be a twisted flag $G$-variety 
for a semisimple group $G$ of inner type over $F$.
Let $p_s$ be a non-trivial rational projector in $\CH^n(X_s\times X_s)$,
i.e., $p_s\circ p_s=p_s$. Then there exists a non-trivial projector $p$
on $\CH^n(X\times X)$ such that $p\times_F F_s=p_s$.
Hence, existence of a non-trivial rational projector $p_s$ on
$\CH^n(X_s\times X_s)$ gives rise to
the decomposition of the Chow motive of $X$
\begin{equation}\label{motproj}
\M(X)\cong(X,p)\oplus (X,\Delta_X-p)
\end{equation}
\end{ntt}

\begin{ntt}[Tautological and quotient bundles] 
Let $A$ be a central simple algebra of degree $n$ over $F$.
Consider a generalized Severi-Brauer variety $\SB_d(A)$ of ideals
of reduced dimension $d$ of $A$. Over the separable closure
it becomes isomorphic to the Grassmannian $\Gr_d(n)$ of $d$-dimensional
planes in a $n$-dimensional affine space.

There is a 
tautological vector bundle 
over $\SB_d(A)$ of rank $dn$ denoted by $\tau_d^A$ and given
by the ideals of $A$ of reduced dimension $d$
considered as vector spaces over $F$.
Over the separable closure $F_s$ this bundle becomes isomorphic to 
$Hom(\E_n,\tau_d)=\tau_d^{\oplus n}$,
where $\E_n$ denotes the trivial bundle of rank $n$ and
$\tau_d$ the tautological bundle over the Grassmannian $\Gr_d(n)$.   

The universal quotient bundle over $\SB_d(A)$, denoted by $\kappa_d^A$,
is the quotient $\E_A/\tau_d^A$ of a trivial bundle $\E_A$ of rank $n^2$
modulo the tautological bundle $\tau_d^A$. Clearly, $\kappa_d^A$ is of rank
$(n-d)n$. Over $F_s$ this bundle
becomes isomorphic to $Hom(\E_n,\kappa_d)=\kappa_d^{\oplus n}$,
where $\kappa_d$ is the universal quotient bundle over $\Gr_d(n)$. 
\end{ntt}

We will extensively use the following fact

\begin{lem}\label{liftbund}
Let $A$ be a central simple algebra over $F$, $r$ a positive integer and 
$B$ a division
algebra which represents the class of the $r$-th tensor power of $A$
in the Brauer group.
Let $X=\SB(A)\times \SB_d(B^\op)$.
Then the bundle
$T_s=\pr_1^*(\tau_1^{\otimes r})\otimes\pr_2^*(\tau_d)$ over $X_s$ 
is a pull-back of some bundle $T$ over $X$. As a consequence,
any Chern class of $T_s$ is a rational cycle.
\end{lem}

\begin{proof} 
According to \cite[10.2]{Pa94} the image of the restriction map on $K_0$
$$
K_0(\SB(A)\times\SB_d(B^\op)) \to K_0(\pl^{n-1}\times \Gr_d(n))
$$
is generated by classes of bundles
$\ind(A^i\otimes B^{-j})\cdot 
[\pr_1^*(\tau_1^{\otimes i})\otimes\pr_2^*(\tau_d^{\otimes j})]$.
\end{proof}

\begin{rem} Observe that the similar fact holds if one replaces
the tensor power by an
exterior (lambda) power. This is due to the fact
that $[\Lambda^r A]=[A^{\otimes r}]$ in $\Br(F)$ (see \cite[10.A.]{KMRT98}).
\end{rem}

\begin{ntt}[Poincar\'e polynomial]\label{poin} By \cite{Fu97} the Poincar\'e polynomial
of a Chow group of a Grassmannian $\Gr_d(n)$ 
is given by the Gaussian polynomial
$$
P(\Gr_d(n),t)={n\choose d}(t)=
\frac{(1-t^n)(1-t^{n-1})\ldots (1-t^{n-d+1})}{(1-t)(1-t^2)\ldots (1-t^d)}.
$$
Observe that for a projective space, i.e., for $d=1$, the respective
polynomial takes the most simple form
$$
P(\pl^{n-1},t)=\frac{1-t^n}{1-t}.
$$
Observe also that if $n$ is a prime integer, 
then the polynomial $P(\pl^{n-1},t)$ always
divides $P(\Gr_d(n),t)$.
\end{ntt}

\section{The subgroup of rational cycles}\label{subrat}

The goal of the present section is to provide an explicit set
of generators for the image of the restriction map 
$\CH(X)\to \CH(X_s)$ modulo $n$, where $X$ is 
a product of a Severi-Brauer variety
by a generalized Severi-Brauer variety corresponding
to division algebras of a prime degree $n$.

\begin{ntt}[Grassmann bundle structure] 
Let $A$ be a division algebra of degree $n$ over $F$, 
$r$ a positive integer
and $B$ a division algebra which represents
the class of the $r$-th tensor power of $A$ in $\Br(F)$.
According to \cite[Proposition~4.3]{IK00} the product 
$\SB(A)\times \SB_d((A^{\otimes r})^\op)$ can be identified
with the Grassmann bundle $\Gr_d(\mathcal{V})$ over $\SB(A)$,
where $\mathcal{V}=(\tau_1^A)^{\otimes r}$ is a locally free sheaf
of (right) $A^{\otimes r}$-modules. 
By Morita equivalence we may replace $A^{\otimes r}$ by $B$ and, hence,
obtain that the product $X=\SB(A)\times \SB_d(B^\op)$
is the Grassmann bundle $\Gr_d(\mathcal{W})$ over $\SB(A)$
for a locally free sheaf of (right) $B$-modules $\mathcal{W}$.
The tautological bundle $T$ over $\Gr_d(\mathcal{W})$
is the bundle
$$
T=\pr_1^*(\mathcal{W})\otimes_B \pr_2^*(\tau_B^\op).
$$

Let $Q=\pr_1^*(\mathcal{W})/T$ denote the universal quotient bundle
over $\Gr_d(\mathcal{W})$.
Over the separable closure it can be identified with
$$
Q_s=\pr_1^*(\tau_1^{\otimes r}) \otimes \pr_2^*(\kappa_d)
$$ 
where $\tau_1$ and $\kappa_d$ are the tautological and quotient
bundles over $\pl^{n-1}$ and $\Gr_d(n)$ respectively.
We shall write this bundle simply as $\tau_1^{\otimes r}\otimes \kappa_d$
meaning the respective pull-backs.
\end{ntt}

\begin{ntt}[Grassmann bundle theorem]
According to Grassmann bundle theorem the Chow ring 
$\CH(X)$ is a free $\CH(\SB(A))$-module
with the basis $\Delta_\lambda(c(Q))$, where $\lambda$ runs through the
set of all partitions $\lambda=(\lambda_1,\ldots,\lambda_d)$ with
$n-d\ge \lambda_1\ge \ldots \ge \lambda_d \ge 0$, $c(Q)$ is the total
Chern class of the quotient bundle $Q$ and $\Delta_\lambda$ is the Schur
function. In other words, for any $k$ 
we have the decomposition
\begin{equation}\label{grbn}
\CH^k(X)\cong 
\bigoplus_{\lambda} \Delta_\lambda(c(Q)) \cdot 
\pr_1^*(\CH^{k-|\lambda|}(\SB(A)))
\end{equation}
\end{ntt}

\begin{ntt}\label{findbasis}
Observe that the decomposition (\ref{grbn}) is compatible with 
the scalar extension $F_s/F$. Hence, any rational cycle $\alpha \in \CH^k(X_s)$
can be represented uniquely as the sum of cycles
\begin{equation}\label{decom}
\alpha = \sum_\lambda \Delta_\lambda(c(Q_s))\cdot \alpha_{\lambda},
\end{equation}
where $\alpha_{\lambda} \in \CH^{k-|\lambda|}(\pl^{n-1})$ is rational.
If $n$ is a prime integer, then all rational cycles of positive codimensions
in $\CH(\pl^{n-1})$
are divisible by $n$ (see \cite[Corollary~4]{Ka95}). 
Hence, considering (\ref{decom}) modulo $n$ we obtain
\end{ntt}

\begin{lem}\label{basis} If $n$ is a prime integer, the cycles
$\Delta_\lambda(c(Q_s))$, where $\lambda$ runs
through the set of all partitions, generate
the subgroup of rational cycles of $\CH(X_s)$ modulo $n$.
In particular case $d=1$ the basis of the subgroup of rational cycles of
$\CH^k(\pl^{n-1}\times \pl^{n-1})$ modulo $n$ consist of
Chern classes $c_k(\tau_1^{\otimes r}\otimes \kappa_1)$, where $k=0\ldots n-1$.
\end{lem}

\section{Applications to Severi-Brauer varieties}\label{sevbrauer}

The following two lemmas were proven by N.~Karpenko 
(see \cite[Theorem~2.2.1]{Ka96} and 
\cite[Criterion~7.1]{Ka00}).
In the present section we provide short proofs of these results restricting
to algebras of prime degrees.

\begin{lem}\label{indecp} The motive of a Severi-Brauer variety of a division algebra
$A$ of a prime degree $n$ is indecomposable.
\end{lem}

\begin{proof}
Consider the product $X=\SB(A)\times\SB(A)$. It can be identified
with the product $\SB(A)\times \SB(B^\op)$, where $[B]=[A^{n-1}]$.
Apply Lemma~\ref{basis} to the case $d=1$ and $r=n-1$.
We obtain that in codimension $k=n-1$ there is
only one basis element of the subgroup of rational cycles modulo $n$
$$
\Delta_{n-1}=c_{n-1}(\tau_1^{\otimes n-1}\otimes \kappa_1)=
\sum_{i=0}^{n-1} (1-n)^{i} \cdot H^i\times H^{n-1-i}\in \CH^{n-1}(X_s)
$$
which is congruent modulo $n$ to the diagonal cycle $\Delta=\sum_{i=0}^{n-1} H^i\times H^{n-1-i}$.

If the motive of $\SB(A)$ splits, then there must exist
a rational projector $p\in \CH^{n-1}(X_s)$ such that $\Delta\pm p$ and $p$
are non-trivial.
By composition rule (\ref{composit}) 
if $p$ is a projector, then
$p=\sum_i \pm H^i\times H^{n-1-i}$, where the index $i$ runs through 
a subset of $\{0\ldots n-1\}$.
From the other hand, since $p$ is rational, 
it must be a multiple of $\Delta$ modulo $n$. 
So the only possibility for $p$ is to coincide either with $\pm\Delta$ or
with $0$, contradiction.
\end{proof}

\begin{lem}\label{seviso} Let $A$ and $B$ be two division algebras of a prime degree $n$
generating the same subgroup in the Brauer group.
Then the motives of $\SB(A)$ and $\SB(B^\op)$ are isomorphic
iff $A=B$ or $B^\op$.
\end{lem}

\begin{proof} Take $r$ such that $[B]=[A^{\otimes r}]$ in $\Br(F)$.
and apply Lemma~\ref{basis} for $X=\SB(A)\times\SB(B^\op)$.
We obtain that the subgroup of rational cycles of codimension
$k=n-1$
is generated by the cycle
$$
\Delta_{n-1}=c_{n-1}(\tau_1^{\otimes r}\otimes\kappa_1)=
\sum_{i=0}^{n-1} (-r)^i \cdot H^i \times H^{n-1-i}\in \CH^{n-1}(X_s).
$$
According to composition rule~(\ref{composit}) and Rost nilpotence
theorem (see \ref{exproj}) any
motivic isomorphisms between $\SB(A)$ and $\SB(B^\op)$ is given
by a lifting of a rational cycle of the kind 
$\sum_{i=0}^{n-1} \pm H^i\times H^{n-1-i} \in \CH^{n-1}(X_s)$ and vice versa.
To finish the proof observe that a rational cycle of this kind
is a multiple of $\Delta_{n-1}$ modulo $n$ iff $r\equiv \pm 1\mod n$.
\end{proof}

\section{Generalized Severi-Brauer varieties}
\label{gsevbrauer}

In the present section we prove the main theorem of the paper
which is formulated as follows

\begin{thm} Let $A$ and $B$ be two division algebras of a prime degree $n$
generating the same subgroup in the Brauer group.
Take an integer $r$ such that $[B]=[A^{\otimes r}]$.
Then the motive of a Severi-Brauer variety $\SB(A)$ is 
a direct summand of the motive of a 
generalized Severi-Brauer variety $\SB_d(B)$
if and only if
$$
d\cdot r \equiv \pm 1 \mod n.
$$
\end{thm}

The cases $d=1$ and $d=n-1$ were considered in Lemma~\ref{seviso}.
From now on we assume $1<d<n-1$. Moreover, by duality we may assume
$d\le \left[\tfrac{n}{2}\right]$.

\paragraph{Proof ($\Rightarrow$)} 
Consider the product $X=\SB(A)\times\SB_d(B^\op)$.
According to Theorem~\ref{basis} the subgroup of rational
cycles of $\CH^k(X_s)$ is generated (modulo prime $n$) by the cycles
$$
\Delta_\lambda =\Delta_\lambda (c(\tau_1^{\otimes r}\otimes\kappa_d)),
$$
where $\lambda$ runs through all partitions with $|\lambda|=k$.
By \cite[Example~A.9.1]{Fu98} we have
$$
\Delta_\lambda = \sum_{\mu\subset\lambda}
d_{\tilde{\lambda},\tilde{\mu}}\cdot c_1(\tau_1^{\otimes r})^{k-|\mu|}\cdot 
\Delta_\mu(c(\kappa_d))=
$$
\begin{equation}\label{cyclg}
=\sum_{i=0}^k (-r)^{k-i}\cdot H^{k-i}\times 
\left( \sum_{\mu\subset\lambda\atop |\mu|=i} 
d_{\tilde{\lambda},\tilde{\mu}}\omega_\mu \right),
\end{equation}
where $\tilde{\lambda}$ denotes the conjugate partition for $\lambda$
, i.e., obtained by interchanging rows and columns in the respective
Young diagram,
$H$ is the class of a hyperplane section of $\pl^{n-1}$, 
$\omega_\mu$ denotes the additive generator of $\CH^{|\mu|}(\Gr_d(n))$ 
corresponding to a partition $\mu$
and the coefficients $d_{\tilde{\lambda},\tilde{\mu}}$ are the binomial determinants
$$
d_{\tilde{\lambda},\tilde{\mu}}=
\left|\left({\tilde{\lambda}_i+n-d-i  \atop \tilde{\mu}_j+n-d-j}\right)\right|_{1\le i,j\le n-d}
$$

Let $k=N$, where $N=d(n-d)$ is the dimension of $\Gr_d(n)$. 
In this codimension
there is only one partition $\lambda$ with $|\lambda|=N$, namely,
the maximal one $\lambda=(n-d,\ldots,n-d)$. 
Let $g$ denote the cycle (\ref{cyclg}) corresponding to this
maximal partition, i.e.,
\begin{equation}\label{hcodim}
g = \sum_{m=0}^{n-1} (-r)^m\cdot H^m\times 
\left( \sum_{|\mu|=N-m} 
d_{\tilde{\lambda},\tilde{\mu}}\omega_\mu \right),
\end{equation}
where $\tilde{\lambda}=(d,d,\ldots,d)$ and the coefficients 
$d_{\tilde{\lambda},\tilde{\mu}}$ 
are given by
$$
d_{\tilde{\lambda},\tilde{\mu}}=
\left|{n-i\choose \tilde{\mu}_j+n-d-j}\right|_{1\le i,j\le n-d}
$$
From now on we denote the coefficient $d_{\tilde{\lambda},\tilde{\mu}}$
by $d_{\mu'}$, where $\mu'$ is the dual partition 
$(n-d-\mu_{d},\ldots,n-d-\mu_1)$. Observe that $|\mu'|=N-|\mu|$.

For an integer $m$ denote by 
$g^{(m)}$ the summand of (\ref{hcodim})
for the chosen index $m$ and by $d_{\mu'}^{(m)}$ the respective
coefficients, i.e.,
$$
g^{(m)}=(-r)^m\cdot 
H^m \times \left(\sum_{|\mu'|=m} d_{\mu'}^{(m)}\omega_\mu\right)
$$

Consider the summands $g^{(0)}$ and $g^{(1)}$.
Since the Chow group $\CH(\Gr_d(n))$ has only one additive generator 
in the last two codimensions $N$ and $N-1$ denoted by $\omega_N$ and $\omega_{N-1}$ respectively, we obtain that
\begin{equation}\label{hcycles}
g^{(0)}=1\times \omega_N\quad\text{and}\quad
g^{(1)}=-rd\cdot H\times \omega_{N-1}
\end{equation}

Now we are ready to finish the {\bf($\Rightarrow$)} part of the theorem.

Assume that the motive of $\SB_d(B^\op)$ contains the motive of $\SB(A)$
as a direct summand. Then 
there must exist two rational cycles 
$$
\alpha\in \CH^{n-1}(\pl^{n-1}\times\Gr_d(n))\quad\text{and}\quad \beta\in 
\CH^N(\Gr_d(n)\times\pl^{n-1})
$$ 
such that
$\beta\circ \alpha=\id$. According to the composition rule (\ref{composit})
this can happen only if the
coefficients before the monomials 
$\omega_N \times 1$ and $\omega_{N-1} \times H$ 
of the cycle $\beta$ are equal to $\pm 1$.
But all rational cycles in codimension $N$ 
are generated (modulo $n$) 
by the transposed cycle $g^t$ which has coefficients
$1$ and $-rd$  before the respective monomials (see (\ref{hcycles})).
This can only be possible if $rd\equiv \pm 1 \mod n$.

\paragraph{Proof ($\Leftarrow$)} 
Assume that the congruence $rd \equiv \pm 1 \mod n$
holds. We want to produce two rational cycles 
$$
\alpha\in \CH^{n-1}(\pl^{n-1}\times\Gr_d(n))\quad\text{and}\quad \beta\in \CH^N(\Gr_d(n)\times\pl^{n-1})
$$ 
such that
$\beta\circ \alpha \in \CH^{n-1}(\pl^{n-1}\times\pl^{n-1})$ is the identity morphism, i.e., the diagonal cycle $\sum_{m=0}^{n-1} H^{n-1-m}\times H^m$. 
This will show that $\M(\SB_d(B^\op))$ contains
$\M(\SB(A))$ as a direct summand.

Assume $rd\equiv 1 \mod n$.
Consider the bundle $\kappa_1\otimes \Lambda^d(\tau_d)$ of rank $n-1$ on the
product $\pl^{n-1}\times\Gr_d(n)$ and define the cycle $f$ as
$$
f=c_{n-1}(\kappa_1\otimes \Lambda^d (\tau_d))=
\sum_{m=0}^{n-1} c_{n-1-m}(\kappa_1)c_1(\tau_d)^m=
$$
$$
=\sum_{m=0}^{n-1} (-1)^m\cdot H^{n-1-m}\times \omega_1^m= \sum_{m=0}^{n-1}
(-1)^m \cdot H^{n-1-m}\times\left( \sum_{|\rho|=m} c_\rho^{(m)} \omega_\rho\right).
$$ 
Observe that $f$ is a rational cycle, since $[A]=[B^{\otimes d}]$.
If $rd\equiv -1 \mod n$, then we take the bundle $\Lambda^{n-d}(\kappa_d)$
instead of $\Lambda^d(\tau_d)$
and obtain the same formulae but without the coefficient $(-1)^m$.

The coefficients $c_\rho^{(m)}$ appearing in the presentation
of $\omega_1^m$ in terms of additive generators $\omega_\rho$
of $\CH^m(\Gr_d(n))$ have the following nice property

\begin{lem}\label{coprime} 
For any partition $\rho$ with $m=|\rho|<n$ the coefficient $c_\rho^{(m)}$
is coprime with $n$.
\end{lem}

\begin{proof}
According to \cite[Example~14.7.11.(ii)]{Fu98} we have
\begin{equation}\label{degreefo}
c_\rho^{(m)}=\deg(a_0,\ldots,a_{d-1})=\frac{m!}{a_0!\,a_1!\ldots a_{d-1}!} \prod_{i>j}(a_i-a_j),
\end{equation}
where the set $a=(a_0,a_1,\ldots,a_{d-1})$ is defined by $a_i=\rho_{d-i}+i$,
$i=0,\ldots,d-1$. Observe that the set $a$ has the property
$0\le a_0<\ldots <a_{d-1}\le n-1$ and corresponds to a Schubert variety
of dimension $m$ which is dual to the Schubert variety of codimension $m$
corresponding to the partition $\rho$ (see \cite[14.7]{Fu98}).
\end{proof}

\begin{rem} The constructed cycle $f\in \CH^{n-1}(\pl^{n-1}\times\Gr_d(n))$ has the
following geometric interpretation.
Consider the Pl\"ucker embedding 
$\SB_d(B) \to \SB(\Lambda^d(B))$. Its graph
defines a correspondence $\Gamma\in \CH(\SB(\Lambda^d(B))\times \SB_d(B))$.
It is known that the motive of $\SB(\Lambda^d(B))$ splits as a direct sum
of twisted motives of $\SB(A)$ (see \cite[Corollary~1.3.2]{Ka96}). Let 
$i:\M(SB(A)) \to \M(\SB(\Lambda^d(B)))$ be the respective splitting.
Then over the separable closure the composition $\Gamma\circ i$ coincides 
with the cycle $f$ corresponding to the case $rd\equiv -1 \mod n$.
After replacing $B$ by $B^\op$ the respective composition $\Gamma\circ i$
will give the cycle $f$ corresponding to the case $rd\equiv 1\mod n$. 
\end{rem}

Consider the transposed cycle $g^t$ introduced in the first part of the proof
$$
g^t = \sum_{m=0}^{n-1} (-1)^m\cdot 
\left( \sum_{|\mu'|=m} 
(r^m d_{\mu'})\cdot \omega_\mu \right)\times H^m
$$

Consider the $m$-th summand of the composition $g^t \circ f$
$$
(g^t\circ f)^{(m)}= 
\left(r^m \sum_{|\mu'|=m} c_{\mu'}^{(m)}\cdot d_{\mu'}^{(m)}\right) \cdot H^{n-1-m}\times H^m
$$
Note that if $rd\equiv -1 \mod n$, 
then the coefficient $(-1)^m$ will appear on the right hand side. 

Assume that the following formulae holds
\begin{equation}\label{equality}
\sum_{|\mu'|=m} c_{\mu'}^{(m)}\cdot d_{\mu'}^{(m)}\equiv d^m \mod n.
\end{equation}
Then the coefficient of $(g^t\circ f)^{(m)}$
is congruent to $1$ modulo  $n$.
We claim that it is possible to modify the cycles $g^t$ and $f$
by adding cycles divisible by $n$ in such a way 
that the coefficient of $(g^t\circ f)^{(m)}$ becomes equal to $1$ 
for each $m$. 

First, we modify the cycle $f$. For each $m$, $0\le m\le n-1$,
we do the following procedure.
Consider the $m$-th summand 
$$
f^{(m)}=(-1)^m \sum_{|\mu'|=m} c_{\mu'} \cdot H^{n-1-m}\times \omega_{\mu'}
$$ 
and its coefficients $c_{\mu'}^{(m)}$.
For $m=0$ and $1$ there is only one additive generator
of $\CH^m(\Gr_d(n))$ ($\omega_0$ and $\omega_1$)
and the respective coefficients are
$c^{(0)}=c^{(1)}=1$. So we set $\alpha^{(0)}=f^{(0)}$ and 
$\alpha^{(1)}=f^{(1)}$. 
For $1<m\le n-1$ the number of generators of $\CH^m(\Gr_d(n))$
is greater than $1$ and all the coefficients 
$c_{\mu'}^{(m)}$ are coprime with $n$, in particular,
they are all non-zero. 
In this case we can modify each $c_{\mu'}^{(m)}$ modulo $n$
by adding cycles of the kind $a\cdot n\omega_{\mu'}$, $a\in \zz$, 
to the cycle $f^{(m)}$ in such a way that 
the greatest common divisor of resulting coefficients, denoted
by $\alpha_{\mu'}^{(m)}$, becomes equal to $1$.
As a result, we obtain a new cycle $\alpha^{(m)}$
having the coefficients $\alpha_{\mu'}$ instead of $c_{\mu'}$.

\begin{dfn}\label{condalpha}
Define a cycle $\alpha$ as $\alpha=\sum_{m=0}^{n-1}\alpha^{(m)}$. 
By construction of $\alpha^{(m)}$ we have
\begin{itemize}
\item $\alpha^{(0)}=\alpha^{(1)}=1$; 
\item $\alpha$ is rational (congruent modulo $n$ to $f$);
\item all the coefficients $\alpha_{\mu'}$ are coprime with $n$;
\item for each $m$ the g.c.d of coefficients $\alpha_{\mu'}^{(m)}$ is $1$;
\item for each $m$ the coefficient of $(g^t\circ \alpha)^{(m)}$ is 
congruent to $1$ modulo $n$.
\end{itemize}
\end{dfn}

Next we modify the second cycle $g^t$. For each $m$ we apply
the following obvious observation
\begin{lem} Let $a_1,\ldots,a_l$ be a finite 
set of integers with $g.c.d.=1$. 
Assume $\sum_i a_i b_i \equiv 1 \mod n$ for some integers $b_i$.
Then there exist integers $b_i' \equiv b_i \mod n$ such that
$\sum_i a_i b_i'=1$.
\end{lem} 
to the congruence
$$
\sum_{|\mu'|=m} c_{\mu'}^{(m)}(r^m d_{\mu'}^{(m)}) \equiv 1 \mod n,
$$
 where $a_i=\alpha^{(m)}_{\mu'}$ are coefficients
of the cycle $\alpha^{(m)}$ and
$b_i=r^m d_{\mu'}^{(m)}$ are coefficients 
of the cycle $(g^t)^{(m)}$. As a result, we obtain a new cycle $\beta^{(m)}$ 
having coefficients $b_i'$ instead of $b_i$.

\begin{dfn}\label{condbeta}
Define a cycle $\beta$ as $\beta=\sum_{m=0}^{n-1}\beta^{(m)}$. By construction 
of $\beta^{(m)}$ we have 
\begin{itemize}
\item $\beta$ is rational (congruent modulo $n$ to $g^t$)
\item for each  $m$ 
the coefficient of $(\beta\circ \alpha)^{(m)}$ is equal to $1$.
\end{itemize}
\end{dfn}

To finish the proof observe that the constructed cycles
$\alpha$ and $\beta$ are rational and the composition
$\beta\circ\alpha$ is the diagonal cycle. 
This implies that the composition $\alpha\circ \beta$ 
is a rational projector which gives rise to the decomposition 
of motives with integral coefficients
(see (\ref{motproj}))
$$
\M(\SB_d(B^\op))\simeq\M(\SB(A))\oplus H
$$
for some motive $H$.

\section{The proof of (\ref{equality})}\label{proofid}

In the present section we prove the following 

\begin{lem}
$$
\sum_{|\mu'|=m} c_{\mu'}^{(m)}\cdot d_{\mu'}^{(m)}\equiv d^m \mod n,
$$
\end{lem}

\begin{proof} First, we express the coefficients $c_{\mu'}$ and
$d_{\mu'}$ in terms of binomial determinants.
According to \cite[Example~14.7.11]{Fu98} we have 
\begin{equation}\label{cpart}
c_{\mu'}=\deg(\omega_1^m \cdot \omega_\mu) = 
\frac{m! \cdot b_0!\, b_1! \ldots b_{d-1}!}{a_0! \,a_1! \ldots a_{d-1} !}\cdot
\left|{a_i \choose b_j}\right|_{0\le i,j\le d-1}
\end{equation}
where $a_i=n-d+i-\mu_{i+1}$ and $b_j=j$ for $i,j=0\ldots d-1$. 
Observe that the sets $a=(a_0,\ldots,a_{d-1})$ and
$b=(b_0,\ldots,b_{d-1})$ with $0\le a_0<\ldots<a_{d-1}\le n-1$ and 
$0\le b_0<\ldots<b_{d-1}\le n-1$ correspond to the classes of 
Schubert varieties $\omega_\mu$ and $\omega_N=\{pt\}$ 
of dimensions $\sum_{i=0}^{d-1} (a_i-i)=m$ and
$\sum_{j=0}^{d-1} (b_j-j)=0$ respectively.

From another hand side, we have 
$$
d_{\mu'}=\left| {n-i \choose \tilde{\mu}_j+n-d-j}\right|_{1\le i.j\le n-d}=
\left|{\tilde{a}_i \choose \tilde{b}_j} \right|_{0\le i,j,\le n-d-1}
$$
where $\tilde{a}_i=d+i$ and $\tilde{b}_j=\tilde{\mu}_{n-d-j}+j$ for $i,j=0\ldots n-d-1$.
Here the sets $\tilde{a}=(\tilde{a}_0,\ldots,\tilde{a}_{n-d-1})$ and 
$\tilde{b}=(\tilde{b}_0,\ldots,\tilde{b}_{n-d-1})$
with $0\le \tilde{a}_0<\ldots \tilde{a}_{n-d-1}\le n-1$ and 
$0\le \tilde{b}_0<\ldots<\tilde{b}_{n-d-1}\le n-1$
correspond to the classes of Schubert varieties on the dual Grassmannian
$\Gr_{n-d}(n)$ of dimensions $N$ and $N-m$ respectively.
Observe that the Schubert variety corresponding to the set $\tilde{b}$ is
dual to the Schubert variety corresponding to the partition $\tilde{\mu}$.
Hence, we have
$$
\omega_1^m = \frac{m! \cdot \tilde{b}_0!\,\tilde{b}_1!\ldots \tilde{b}_{n-d-1}!}{\tilde{a}_0!\,\tilde{a}_1!\ldots \tilde{a}_{n-d-1}!}d_{\mu'}\cdot \omega_{\tilde{b}}
$$
And by duality we obtain
$$
\frac{m! \cdot \tilde{b}_0!\,\tilde{b}_1!\ldots \tilde{b}_{n-d-1}!}{\tilde{a}_0!\,\tilde{a}_1!\ldots \tilde{a}_{n-d-1}!}d_{\mu'}=
\frac{m! \cdot b_0!\, b_1! \ldots b_{d-1}!}{a_0! \,a_1! \ldots a_{d-1} !}
\cdot \left|{a_i \choose b_j}\right|_{0\le i,j\le d-1}
$$
Since the integers $a_0,\ldots,a_{d-1}$ form the complement
of $\tilde{b}_0,\ldots,\tilde{b}_{n-d-1}$ in the set of integers from $0$ to $n-1$ (see \cite[Example~14.7.5]{Fu98}), we obtain
\begin{equation}\label{dpart}d_{\mu'}=\left|{a_i \choose b_j}\right|_{0\le i,j\le d-1}
\end{equation}

According to \cite[Example~14.7.11]{Fu98} we express the binomial
determinant appearing in (\ref{cpart}) and (\ref{dpart}) 
in terms of Vandermonde determinant 
$D_a=\prod_{i<j}(a_j-a_i)$ 
and get
\begin{equation}\label{cdexp}
c_{\mu'}=\frac{m!\cdot D_a}{a_0!\,a_1!\ldots a_{d-1}!}, \qquad\qquad
d_{\mu'}=\frac{D_a}{0!\,1!\ldots (d-1)!}.
\end{equation}

Then the formulae (\ref{equality}) we want to prove 
turns into 
\begin{equation}\label{cong}
\frac{m!}{0!\,1!\ldots (d-1)!}\cdot \sum_a  \frac{D_a^2}{a_0! a_1!\ldots a_{d-1}!} \equiv d^m \mod n,
\end{equation}
where the sum is taken over all sets of integers
$a=(a_0,\ldots,a_{d-1})$ such that $0\le a_0<\ldots<a_{d-1}\le n-1$ and
$\sum_{i=0}^{d-1} a_i=m+\tfrac{d(d-1)}{2}$.

\paragraph{The case $d=2$.}
In this case (\ref{cong}) follows from the following elementary fact
\begin{equation}\label{cong2}
\frac{m!}{2} \cdot 
\sum_{0\le x_1,x_2 \atop x_1+x_2=m+1} \frac{(x_1-x_2)^2}{x_1! x_2!}=2^m
\end{equation}
To prove it consider the following chain of identities
$$
\frac{m! \cdot (x_1-x_2)^2}{2\cdot x_1! x_2!}=
\frac{(m+1)!}{x_1! (m+1-x_1)!}\cdot 
\left(\frac{m+1}{2}-\frac{2x_1(m+1-x_1)}{m+1} \right)=
$$
$$
=\frac{m+1}{2}{m+1 \choose x_1}-2m{m-1 \choose x_1-1}
$$
Then taking the sum we obtain the desired equality
$$
\sum_{x_1=0}^{m+1} \frac{m+1}{2}{m+1 \choose x_1}-2m{m-1 \choose x_1-1} =
\frac{m+1}{2}\cdot 2^{m+1} - 2m\cdot 2^{m-1}=2^m
$$

Observe that for $m<n-1$ the left hand side of (\ref{cong}) coincides
with the left hand side of (\ref{cong2}), hence,
we obtain the equality in (\ref{cong}) (not just a congruence modulo $n$). 
For $m=n-1$ the left hand side of (\ref{cong}) is equal to
$2^{n-1}-n$.

\paragraph{The general case.}
It turns out that the identity (\ref{cong2}) is a particular case
of the following combinatorial identity known as Robinson-Schensted correspondence
(see \cite[4.3.(5)]{Fu97})
\begin{equation}\label{RS}
\sum_{|\xi|=m} d_{\xi}(d)\cdot f^{\xi} = d^m,
\end{equation}
where the sum is taken over all partitions 
$\xi=(\xi_1 \ge \xi_2 \ge \ldots \ge \xi_d\ge 0)$  with
$|\xi|=m$,
$d_\xi(d)$ denote the number of Young tableaux on the shape $\xi$
whose entries are taken from the set $(1,\ldots,d)$ and
$f^\xi$ denote the number of standard tableaux on the shape $\xi$.

By using Hook length formulae (see \cite[4.3, Exercise~9]{Fu97})
we obtain
\begin{equation}\label{cl}
f^\xi=\frac{m!\cdot D_l}{l_0 !\ldots l_{d-1} !},
\end{equation}
where $l=(l_0,\ldots,l_{d-1})$ is a strictly increasing set of non-negative
integers defined from the partition $\xi=(\xi_1,\ldots,\xi_d)$
by $l_{d-i}=\xi_i+d-i$ and $D_l$ is the Vandermonde determinant for $l$. 
By definition we have $\sum_{i=0}^{d-1} (l_i-i)=m$.

By \cite[Corollary~13]{GV85} we have
\begin{equation}\label{dl}
d_\xi(d)=
\left|{l_i \choose j} \right|_{0\le i,j,\le d-1}=
\frac{D_l}{0!\,1!\ldots (d-1)!}
\end{equation}

Observe that if the set $l$ is bounded by $n-1$, i.e.,
$l_{d-1}\le n-1$, the expressions (\ref{cl}) and $(\ref{dl})$ coincide
with the expressions (\ref{cdexp}) defining
the coefficients $c_{\mu'}$ and $d_{\mu'}$ 
respectively (take $\xi=\mu'$ and $a=l$). 

Assume that $l_{d-1} \ge n$. 
Then $l_{d-2}\le n-1$, i.e., $l_{d-1}$ is the only element of $l$ which
exceeds $n-1$. Indeed, if this is not the case then we have a sequence
of inequalities
$$
n-1+\tfrac{d(d-1)}{2}\ge m+\tfrac{d(d-1)}{2}=0+1+\ldots+l_{d-3}+l_{d-2}+l_{d-1}\ge 
$$
$$
\ge 0+1+\ldots + (d-3) +  n + (n+1)=\tfrac{(d-3)(d-2)}{2}+2n+1
$$
which can be rewritten as
$$
d\ge \frac{n+5}{2}.
$$
But we have assumed from the beginning that $d\le \left[\tfrac{n}{2}\right]$,
contradiction.

Moreover, by the similar arguments 
one can check that $l_{d-1}<2n$. 

Now consider the product $x=d_\xi(d)\cdot f^\xi$, when $\xi$ is a partition
for which the respective $l_{d-1}\ge n$.
Since $n$ is prime, the denominator of $x$ is divisible by $n$ 
but not by $n^2$. Since $x$ is an integer, the numerator of $x$ must
be divisible by $n$ as well. From this we conclude that 
$D_l$ is divisible by $n$. But the numerator of $x$ is the product
of the square of $D_l$ by something, hence, it is divisible by $n^2$.
So we obtain that $x$ must be divisible by $n$.  
This means that modulo $n$ the left hand side of (\ref{RS}) is congruent
to the left hand side of (\ref{cong}).
\end{proof}

\section{Grassmannian $\Gr(2,n)$}\label{scop}

The goal of the present 
is to extend Theorem~\ref{introthm} to the case
of algebras of an arbitrary odd degree $n\ge 5$ and $d=2$

\begin{thm}\label{gr2}
Let $A$ and $B$ be two central simple algebras of an odd degree $n$
over a field $F$ generating the same subgroup in the Brauer group $\Br(F)$.
Then the motive of a Severi-Brauer variety $\SB(A)$ is 
a direct summand of the motive of a generalized Severi-Brauer variety 
$\SB_2(B)$
if and only if
$$
[A]=\pm 2\,[B] \text{ in } \Br(F).
$$
\end{thm}

\paragraph{Proof} 
Consider the cycles $g$ and $f$ defined in Section~\ref{gsevbrauer}.
Clearly, $g$ and $f$ are rational, since Lemma~\ref{liftbund}
holds for any algebras.

\subparagraph{($\Rightarrow$)} 
Repeating the arguments
of \ref{findbasis} for an odd integer $n$ and $d=2$
one obtains (in view of \cite[Corollary~4]{Ka95}) that the subgroup
of rational cycles of $\CH(X_s)$ modulo $n$ is generated
in codimension $N=\dim\Gr_d(n)$ by the cycles
$$
\tfrac{n}{(n,N-|\lambda|)}\cdot 
 H^{N-|\lambda|}\times 1 \cdot \Delta_\lambda(c(Q_s)),
$$
for all partitions $\lambda$ with $N-(n-1)\le |\lambda|\le N$.
Observe that in the case $|\lambda|=N$ one obtains precisely the cycle $g$.

Since the coefficient $\tfrac{n}{(n,N-|\lambda|)}$ is divisible
by $n$ when $|\lambda|=N-1$, we may exclude the cycles
with $|\lambda|= N-1$ from the set of generators.
The latter implies that the last two
non-trivial monomials (see (\ref{hcycles})) of any rational
cycle in $\CH^N(X_s)$ come only from the cycle $g$. 
And we finish the proof as in Section~\ref{gsevbrauer}.

\subparagraph{($\Leftarrow$)}
We claim that it is still possible to modify cycles $f$ and $g^t$
modulo $n$ in such a way that the obtained rational cycles $\alpha$
and $\beta$ will satisfy the property $\beta\circ\alpha=\id$, hence,
providing the desired motivic decomposition for $\M(\SB_2(B))$.

The following lemma allows us to take $\alpha=f$ 
(see Definition~\ref{condalpha}).
\begin{lem} Fix a codimension $m$, $0\le m\le n-1$.
For any partition $\rho$ with $|\rho|=m$ let $c_\rho^{(m)}$ be the coefficient
appearing in Lemma~\ref{coprime}.
Then the greatest common divisor of all the coefficients $c_\rho^{(m)}$
is $1$.
\end{lem}

\begin{proof}
Since $d=2$, for any codimension $m$ which is less than $n-1$,
there is a partition $\rho=(m,0)$ with $c_\rho^{(m)}=1$. 
In the last codimension $m=n-1$ consider
two partitions $\rho=(1,n-2)$ and $\rho'=(2,n-3)$.
By (\ref{degreefo}) 
the respective coefficients $c_{\rho}^{(n-1)}$ and $c_{\rho'}^{(n-1)}$
are equal to $\deg(1,n-1)=n-2$ and $\deg(2,n-2)=\tfrac{(n-1)(n-4)}{2}$.
Since $n$ is odd, $n-2$ and $\tfrac{(n-1)(n-4)}{2}$ are coprime.
\end{proof}

Choose $\beta$ as in the Section~\ref{gsevbrauer}. To finish
the proof we have to prove the congruence 
(\ref{equality}).
But this was done already in Section~\ref{proofid}, where
we treated the case $d=2$. Namely, we proved that for $m<n-1$
there is an exact equality (not just a congruence modulo $n$)
and for $m=n-1$ the left hand side of (\ref{equality}) is,
indeed, equal to $2^{n-1}-n$.
This finishes the proof of Theorem~\ref{gr2}.

\subparagraph{Acknowledgements} I would like to thank
B.~Calm\'es, N.~Karpenko, I.~Panin, M.~Rost and N.~Semenov for stimulating
discussions on the subject of the present paper.
I am grateful to
Max-Planck-Institut f\"ur Mathematik
in Bonn and University of Bielefeld for hospitality and support.
This work was partially supported by RTN-Network HPRN-CT-2002-00287.

\bibliographystyle{chicago}

\end{document}